\documentclass[a4paper,11pt]{amsart}
\usepackage{amsmath, amsthm, amssymb, amsfonts}
\usepackage{amscd}
\usepackage[all]{xy}
\usepackage{graphicx,color}

\newtheoremstyle{Style}%
{.5em}{.5em}%
{\it}%
{}%
{\sc}%
{\ {\bf---}}%
{ }%
{}%

\newtheoremstyle{StyleRemarque}%
{.5em}{.5em}%
{\it}%
{}%
{\slshape}%
{.\ }%
{ }%
{}%

\theoremstyle{Style}
\newtheorem{defn}{Definition}[section]
\newtheorem{lem}[defn]{Lemma}
\newtheorem{prop}[defn]{Proposition}
\newtheorem{thm}[defn]{Theorem}
\newtheorem{cor}[defn]{Corollary}
\newtheorem*{thm*}{Theorem}

\theoremstyle{StyleRemarque}
\newtheorem{rem}{Remark}
\newtheorem*{rem*}{Remark}

\newtheorem*{ex*}{Example}

\newcounter{mnotecount}[section]
\renewcommand{\themnotecount}{\thesection.\arabic{mnotecount}}
\newcommand{\mnote}[1]
{\protect{\stepcounter{mnotecount}}$^{\mbox{\footnotesize  $
      \bullet$\themnotecount}}$ \marginpar{\raggedright\tiny\em
    $\!\!\!\!\!\!\,\bullet$\themnotecount: #1} }

\newcommand{\subup}[3]{\underset{#1}{\overset{#2}{#3}}}
\newcommand{\abs}[1]{\left\vert#1\right\vert}

\newcommand{\set}[1]{\left\{#1\right\}}
\newcommand{\Forall}[2]{\forall \, #1 \in #2, \:}

\newcommand{\diagcommut}[8]{
$$\begin{CD}
 #1@>{#2}>>#3  \\
 @VV{#4}V   @VV{#5}V \\
 #6@>{#7}>>#8 \\
\end{CD}$$}
\newcommand{\restric}[1]{\vert_{#1}}
\newcommand{\appli}[3]{#1 \, : \, #2 \To #3}

\newcommand{\rl}{\mathbb R}
\newcommand{\cx}{\mathbb C}
\newcommand{\ir}{\mathbb Z}

\newcommand{\sph}{\mathbb S}

\newcommand{\To}{\longrightarrow}
\newcommand{\oo}{\mathcal O}
\newcommand{\dd}{\mathcal D}
\newcommand{\ee}{\mathcal E}
\newcommand{\nn}{\mathcal N}
\newcommand{\kk}{\mathcal K}
\renewcommand{\d}{\operatorname{d}}

\DeclareMathOperator{\tr}{Tr}
\DeclareMathOperator{\Rm}{Rm}
\DeclareMathOperator{\ric}{Ric}
\DeclareMathOperator{\scal}{Scal}

\DeclareMathOperator{\vol}{vol}

\DeclareMathOperator{\grad}{grad}
\DeclareMathOperator{\hess}{Hess}

\DeclareMathOperator{\dirac}{\not{\hskip -1mm D}} 
\DeclareMathOperator{\Div}{div}
\DeclareMathOperator{\pre}{Re}
\DeclareMathOperator{\pim}{Im}

\begin{document}

\title[A mass for asymptotically complex hyperbolic manifolds.]{\bfseries{A mass for asymptotically complex hyperbolic manifolds.}}
\date{\today}
\thanks{\textit{2000 Mathematics Subject Classification.}  }
\thanks{The authors benefit from the ANR grant \textit{GeomEinstein}, ANR-06-BLAN-0154-02.}
\email{daniel.maerten@yahoo.fr, minerbe@math.jussieu.fr}
\address{Universit\'e Paris 6, Institut de Math\'ematiques de Jussieu, UMR CNRS 7586, 175 rue du Chevaleret, 75013 Paris, France.}
\maketitle

\begin{abstract}
We prove a positive mass theorem for complete K\"ahler manifolds 
that are asymptotic to the complex hyperbolic space.
\end{abstract}

\tableofcontents

\section*{Introduction}
The aim of this paper is to provide a positive mass theorem in the realm of asymptotically complex hyperbolic K\"ahler manifolds, extending previous results by M. Herzlich \cite{Her} and Boualem-Herzlich \cite{BH}. Before explaining them, let us recall the history of the subject. 

The classical positive mass theorem finds its roots in general relativity \cite{ADM} and deals with asymptotically Euclidean manifolds, namely complete Riemannian manifolds $(M^n,g)$ whose geometry at infinity tends to that of the flat Euclidean space : $M$ is diffeomorphic to $\rl^n$ outside a compact set and $g$ goes to $g_{\rl^n}$ at infinity. Under a nonnegativity assumption on the curvature ($\scal_g \geq 0$), the positive mass theorem roughly asserts that such manifolds possess a global Riemannian invariant, which is called a mass, which is obtained by computing the limit of integrals over larger and larger spheres, which is a non-negative number and which vanishes only when the manifold is isometric to the model flat space. This ``Euclidean mass'' is given in some chart at infinity by 
\begin{equation}\label{AEmass}
\mu_g = - \frac14 \lim_{R \To \infty} \int_{\sph_R}  * \left( \Div g + d\tr g \right),
\end{equation}  
where the sphere, Hodge star, the divergence and the trace are defined with respect to the Euclidean metric $g_{\rl^n}$ at infinity.
It is not obvious that this quantity depends only on $g$ and not on the chart but \cite{Bart} proved it. We refer to \cite{SY1,SY2,Wit,Bart,LP} for details and ``classical'' proofs and to \cite{Loh} for a more recent and more general treatment. The mathematical interest of such a theorem is the rigidity result it involves: under a nonnegativity assumption on the curvature, it asserts a model metric at infinity (the Euclidean metric, here) cannot be approached at any rate, the obstruction being precisely the mass. Another striking feature of this theorem is its role as a key step in the first proof of the Yamabe theorem \cite{LP} ; see also \cite{ST} for a nice application of the positive mass theorem to obtain a rigidity result for compact manifolds with boundary.

It is possible to extend these ideas behind this theorem in several contexts, involving different models at infinity. The authors of \cite{MO,AD,CH} studied the case of manifolds whose model at infinity is the real hyperbolic space $\rl H^n$. More sophisticated models are also considered in \cite{Dai,Min}. The notion of real hyperbolic mass introduced in \cite{CH} is at the very root of our work so we need to explain what it looks like (note it is more general than what follows). Basically, while the Euclidean mass (\ref{AEmass}) is a single number, the real hyperbolic mass appears naturally as a linear functional on some finite dimensional vector space $\nn$, attached to the model at infinity, $\rl H^n$. More precisely, \cite{CH} defines $\nn$ as the set of functions $u$ on $\rl H^n$ such that $\hess_{\rl H^n} u = u g_{\rl H^n}$ and the mass linear functional is given by
\begin{equation}\label{ARHmass}
\mu_g (u) = - \frac14 \lim_{R \To \infty} \int_{\sph_R}  * \left[ \left( \Div g + d\tr g \right) u 
- \tr ( g- g_{\rl H^n} ) du + (g-g_{\rl H^n} ) (\grad u ,.) \right],
\end{equation}  
where the right-hand side is computed with respect to the real hyperbolic metric $g_{\rl H^n}$ at infinity. It turns out that $\nn$ can be interpreted as the set of parallel sections for some natural vector bundle $\ee$, endowed with a natural connection $\nabla^{RH}$ preserving a natural Lorentz metric $h$ : if you see $\rl H^n$ as a hypersurface in Minkowski space $\rl^{n,1}$,
then $(\ee,h)$ is simply the restriction of the tangent bundle $T \rl^{n,1}$ and $\nabla^{RH}$ is induced by the flat connection on Minkowki space ; from this point of view, $\nn$ identifies naturally with $\rl^{n,1}$. We will give more details on this picture in the first section of this text and explain why it is natural that the mass belongs to the dual of this space $\nn$ and why the formula above is a relevant geometric invariant. The basic idea is that \emph{any} Riemannian manifold carries a natural connection $\nabla^{RH}$ that is flat iff the manifold is (locally) hyperbolic ; the mass appears naturally as an obstruction to the construction of parallel sections for this connection.

The model at infinity we wish to consider is the complex hyperbolic space $\cx H^m$, which is a counterpart of the real hyperbolic space $\rl H^n$ in complex geometry. Up to scale, it is indeed the unique simply-connected complete K\"ahler manifold with 
constant \emph{holomorphic} sectional curvature. To be concrete, let us recall that $\cx H^m$ can be seen as the unit ball in $\cx^m$, endowed with its standard complex structure $J$ and with the K\"ahler metric 
$$
g_{\cx H^m} = g_{\cx H^m} = \frac{1}{(1-s^2)^2} \left( ds^2 + (J ds)^2 + s^2(1-s^2) g_{FS}   \right),
$$
where $s$ is the radial coordinate in $\cx^m$ and $g_{FS}$ is the Fubini-Study metric of $\cx P^{m-1}$, pulled-back to $\cx^m \backslash \{0\}$. Setting $s= \tanh r$, we obtain $\cx H^m$ as the complex manifold $\cx^m$ endowed with 
$$
g_{\cx H^m} = dr^2 + (\sinh 2r)^2 \eta^2 + (\sinh r)^2 g_{FS}
$$
where $\eta = -\frac{J dr}{\sinh 2r}$ is the standard contact form on $\sph^{2m-1}$. This is an Einstein metric with scalar curvature $-4m(m+1)$ -- the holomorphic sectional curvature is $-4$. The most useful description of $\cx H^m$ for us is yet another one. Let $\cx^{m,1}$ denote the vector space $\cx^{m+1}$ endowed with a Hermitian form $h$ of (complex) signature $(m,1)$. Then the level set $h=-1$ in $\cx^{m,1}$, endowed with the induced metric, is a well-known Lorentz manifold, called Anti-de-Sitter space $AdS^{2m,1}$.
The quotient of $AdS^{2m,1}$ by the scalar action of $\sph^1$ is then a Riemannian manifold and it is precisely $\cx H^m$.  

In this paper, we define \emph{asymptotically complex hyperbolic manifolds} as complete K\"ahler manifolds $(M^{2m},g,J)$ such that:
\begin{enumerate}
	\item[(i)] $M$ minus a compact subset is biholomorphic to  $\cx H^m$ minus a ball and,
	\item[(ii)] through this identification, $g - g_{\cx H^m} = \oo(e^{-a r})$ with $a > m + \frac12$ (in $C^{1,\alpha}$).
\end{enumerate}
Note the definition in \cite{Her}, while apparently weaker, is indeed equivalent (cf. the remark after Definition \ref{defasympcx}). The papers \cite{Her,BH} prove rigidity results about asymptotically complex hyperbolic manifolds that look like the rigidity part (the vanishing mass part) of a positive mass theorem. What is the mass in this setting ? 

In complete analogy with the real hyperbolic case, we will see the mass as a linear functional on some natural 
finite-dimensional vector space $\nn$ attached to the model at infinity, $\cx H^m$. This vector space $\nn$ is best described as a
set of parallel sections for some natural connection $\nabla^{CH}$ on some natural vector bundle $\ee$ over $\cx H^m$. The vector bundle $\ee$ is indeed $\Lambda^2_J \cx H^m \oplus T \cx H^m \oplus \rl$ and the relevant connection comes from the flat connection on $\cx^{m,1}$. Details will be given in the text. To keep this introduction short, let us just point out that $\nn$ identifies naturally to the vector space $\Lambda^2_{J} \cx^{m,1}$ of $J$-invariant $2$-forms on $\cx^{m,1}$ and also admits a description as the set of functions $u$ satisfying a natural third-order equation. Then we define the complex hyperbolic mass by
\begin{equation}\label{ACHmass}
\mu_g (u) = - \frac14 \lim_{R \To \infty} \int_{\sph_R}  * \left[ \left( \Div g + d\tr g \right) u 
- \frac12 \tr ( g- g_{\cx H^m} ) du \right],
\end{equation}  
where everything on the right-hand side is computed with respect to the complex hyperbolic metric $g_{\cx H^m}$ at infinity.
The following positive mass theorem holds in this context.

\begin{thm}
Let $(M,g,J)$ be a spin asymptotically complex hyperbolic manifold with $\scal_g \geq \scal_{g_{\cx H^m}}$. When 
the complex dimension of $M$ is even, we also assume that $M$ is contractible. Then $\mu_g$
is a well defined linear functional, up to an automorphism of the model. It vanishes if and only if $(M,g,J)$ is the complex hyperbolic space. 
\end{thm}

The mass also satisfies a nonnegativity property : it takes non-negative values on some distinguished orbits of the action of $PU(m,1)$. Under these assumptions, the mass may very well take infinite values : we then decide that is infinite and the finiteness of the mass does not depend on the choice of the chart at infinity. A simple criterion for the mass to be finite is that the rate of fall-off to the model metric at infinity be sufficiently fast : $a \geq 2m +1$. A consequence of the positive mass theorem is that if $a > 2m+1$, then $(M,g,J)$ is complex hyperbolic ; it is the object of \cite{Her,BH}. Our viewpoint yields a different proof, maybe more direct.

The spin assumption is classical : we need spinors to implement Witten's techniques \cite{Wit}. The additional topological assumption in the even dimensional case is quite technical but could be weakened ; we refer to the text for a more precise statement. This assumption was already used in \cite{BH}. Basically, Witten's techniques require some ``special'' spinors to exist on the model.
Complex hyperbolic spaces of odd dimension possess such distinguished spinors, called K\"ahlerian Killing spinors. Such spinors do not exist in the even dimensional case, so we need an extra trick. This has a cost : an additional assumption.

Let us make a general remark. The previous positive mass theorems were very related to -- if not completely immersed into -- physical ideas. The asymptotically complex hyperbolic realm does not share this feature (yet). We feel it is all the more interesting to observe that the very physical idea of a positive mass theorem carries over to purely mathematical settings. It might indicate that there is a general mechanism, waiting for new applications. 

The structure of the paper is as follows. In a first section, we will discuss the real hyperbolic case and explain 
how the introduction of a natural hyperbolic connection gives a slightly different proof of the hyperbolic positive 
mass theorem in \cite{CH}. This case will also serve as a helpful guide for the complex case. In a second section, we 
will introduce  the complex hyperbolic connection and describe its basic features and in particular how it interacts with 
the so-called K\"ahlerian Killing spinors. A third section is devoted to the proof of the complex hyperbolic positive
mass theorem in odd complex dimensions. A short fourth section will describe how to extend these arguments to the even 
dimensions. Finally, an appendix describes an example.

\emph{Acknowledgements.}
The authors would like to thank Marc Herzlich for bringing the problem to their attention, but also Olivier Biquard,
Elisha Falbel and Paul Gauduchon for useful discussions.

\section{The real hyperbolic case.}

In this section, we briefly review some aspects of \cite{CH}. The result mentioned here is not new, but we feel
that the reformulation we propose might be useful. It is indeed quite simple and turns out to generalize to the 
complex hyperbolic setting. The basic question we address is the following : given a Riemannian manifold that looks like the hyperbolic space, how can one ensure that it is actually the hyperbolic space ? The approach suggested here consists in finding a connection characterizing the hyperbolic space, in that it is flat only on a (locally) hyperbolic manifold, and then try to build parallel sections for this connection.

Let us start with a Riemannian manifold $(M^n,g)$. We consider the vector bundle  $\ee = T^*M \oplus \rl$ obtained as the sum
of the cotangent bundle and of the trivial real line bundle. It can be endowed with a natural Lorentz metric $h$ : for $\alpha \in T^*_x M$ and $u\in \rl$, we set 
$$
h(\alpha,u) = |\alpha|^2_g - u^2.
$$
We will say that an element $(\alpha,u)$ of $\ee$ is future light-like if $h(\alpha,u) =0$ and $u>0$. We then define a connection $\nabla^{RH}$ on $\ee$ : if $(\alpha,u)$ is a section of $\ee$ and $X$ a vector field on $M$,
$$
\nabla^{RH}_X 
\begin{pmatrix}
\alpha \\ u 
\end{pmatrix}
:= 
\begin{pmatrix}
\nabla_X^g \alpha - u g(X,.) \\
d_X u - \alpha (X)
\end{pmatrix}.
$$
This connection is metric with respect to $h$. Moreover, an easy computation shows that its curvature vanishes 
if and only if $g$ has sectional curvature $-1$, which is why we call this connection \emph{hyperbolic}. In case $(M,g)$ is $(\rl H^n,g_{\rl H^n})$, this construction is clear : $\rl H^n$ is embedded into Minkowski space $\rl^{n,1}$ as the hypersurface
$$
\set{ x \in \rl^{n+1} \, / \; x_1^2 + \dots + x_n^2 - x^2_{n+1} =-1 \; \text{ and } x_{n+1} >0} ;
$$
then $\ee = \ee_{\rl H^n}$ identifies with $T \rl^{n,1} \restric{\rl H^n}$, $h$ is induced by the Minkowski metric and $\nabla^{RH}$ by the flat Minkowski connection.

\begin{rem}
Note that this construction admits an obvious spherical analogue. On the same vector bundle, one can consider the obvious positive definite metric and change a sign in the formula for the connection to make it metric, which results in a \emph{spherical} connection
: it is flat if and only if $(M,g)$ is locally isometric to the sphere with constant sectional curvature $+1$. 
\end{rem}

Observe also that a section $(\alpha,u)$ of $\ee$ is parallel for this connection
if and only if $\alpha = du$ and $\hess_g u = u g$. We call $\nn$ the space of parallel
sections for $\nabla^{RH}$ and $\nn^+$ the subset of future light-like elements of $\nn$.
Observe that in the model case, $\nn_{\rl H^n}$ identifies with $\rl^{n,1}$ and
$\nn_{\rl H^n}^+$ is simply the future isotropic half-cone.

\begin{rem}
The spherical analogue yields the equation $\hess_g u = - u g$, which is called Obata equation and has been much studied
(\cite{Oba,Gal}), in relation with the bottom of the spectrum of Riemannian manifolds with $\ric \geq \ric_{\sph^n}$.
Obata equation admits a non-trivial solution only on the standard sphere.
\end{rem}

We are interested in trying to produce parallel sections of $\ee$. To do this, we follow an indirect path, assuming $M$ to be spin and looking at a related spinorial connection. An imaginary Killing spinor $\psi$ on a spin Riemannian manifold $(M,g)$ is a section of the spinor bundle  such that
$$
\nabla^g \psi + \frac i2 \psi = 0.
$$
The space of Killing spinors is denoted by $\kk$. The following crucial observation follows from a straightforward computation :
\begin{equation}\label{impliqueRH}
\nabla^g \psi + \frac i2 \psi = 0 \quad \Rightarrow \quad 
\nabla^{RH} \begin{pmatrix}
d \abs{\psi}_g^2 \\ \abs{\psi}_g^2 
\end{pmatrix}
=0.
\end{equation}
In other words, we obtain a map $\appli{Q}{\kk}{\nn}$ by setting $Q(\psi) = (d \abs{\psi}^2 , \abs{\psi}^2)$.

To understand the relevance of imaginary Killing spinors, we must describe them on the model $\rl H^n$ : basically, they are induced by the constant spinors of $\rl^{n,1}$ and they trivialize the spinor bundle of $\rl H^n$ \cite{CH}. Using explicit formulas (cf. \cite{CH}),
one can see that $\nn_{\rl H^n}^+$ lies inside $Q_{\rl H^n}(\kk_{\rl H^n})$ ; indeed, it is sufficient to prove that one of the Killing spinors is mapped into $\nn_{\rl H^n}^+$ and then use the equivariance of $Q$ with respect to the natural actions of $O^+(n,1)$, together with the transitivity of $O^+(n,1)$ on $\nn_{\rl H^n}^+$. 

Now assume the complete Riemannian manifold $(M,g)$ is asymptotically hyperbolic in the following sense : $M$ minus a compact is diffeomorphic to $\rl H^n$ minus a ball and, through this diffeomorphism, $g = g_{\rl H^n} + \oo(e^{-a r})$ with $a > \frac{n+1}{2}$ in $C^{1,\alpha}$ ($r$ is the distance to some point in $\rl H^n$). We further assume that the scalar curvature $\scal_g$ of $(M,g)$ is greater than or equal to 
$\scal_{g_{\rl H^n}}$. Then \cite{CH} proves that for any Killing spinor $\psi$ on $\rl H^n$, one can find a unique spinor $\tilde{\psi}$ that is asymptotic to $\psi$ and harmonic for some natural Dirac operator. A Witten's like argument, based on Lichnerowitz formula, then leads to 
\begin{equation}\label{ippRH}
\int_M \left( \abs{\nabla^g \tilde{\psi} + \frac i2 \tilde{\psi}}^2 + \frac14 \left(  \scal_g - \scal_{g_{\rl H^n}}\right) \abs{\tilde{\psi}}^2 \right) = \lim_{R \to \infty} \int_{\sph_R} \dots
\end{equation}
The right hand side can be computed explicitly : it is  
\begin{eqnarray*}
\mu_g \left( Q_{\rl H^n}(\psi) \right) &=& - \frac14 \lim_{R \To \infty} \int_{\sph_R}  * \Big[ \left( \Div g + d\tr g \right) \abs{\psi}^2 \\
& &- \tr ( g- g_{\rl H^n} ) du + (g-g_{\rl H^n} ) (d\abs{\psi}^2 ,.) \Big],
\end{eqnarray*}
where the sphere, the Hodge star, the divergence, the trace and the identification between vectors and forms on the right-hand side are taken with respect to $g_{\rl H^n}$. Since $\nn_{\rl H^n}^+ \subset Q_{\rl H^n}(\kk_{\rl H^n})$, it follows that $\mu$ yields a linear functional on the linear span of 
$\nn_{\rl H^n}^+$, namely on $\nn_{\rl H^n} \cong \rl^{n,1}$. In view of (\ref{ippRH}), it is clearly non-negative on $\nn_{\rl H^n}^+$. Besides, if $\mu$ vanishes, the left-hand side of (\ref{ippRH}) is zero, so that every $\tilde{\psi}$ is a Killing spinor on $(M,g)$ ; then, in view of (\ref{impliqueRH}), for any element $\sigma=Q_{\rl H^n}(\psi)$ of $\nn_{\rl H^n}^+$, there is an element $\tilde{\sigma} = Q_g(\tilde{\psi})$ of $\nn_g$ that is asymptotic to $\sigma$. As a consequence, $\nn_g$ has maximal dimension, which implies that $(M,g)$ is (locally) hyperbolic. In view of its asymptotic, it is bound to be $\rl H^n$. 

\begin{thm}[\cite{CH}]
Let $(M^n,g)$ be an asymptotically hyperbolic spin manifold, with $\scal_g \geq \scal_{g_{\rl H^n}}$. Then the
linear functional $\mu_g$ on $\nn_{\rl H^n}$ introduced above is well-defined up to an automorphism of $\rl H^n$,
it is non-negative on $\nn_{\rl H^n}^+$ and it vanishes iff $(M^n,g)$ is isometric to $\rl H^n$.
\end{thm}

The fact that the orbit of $\mu_g$ under the action of $O^+(n,1)$ does not depend on the chart at infinity is not obvious but is proved in \cite{CH}. Under our assumptions, $\mu_g$ may take infinite values (i.e. formula (\ref{ippRH}) may be $+\infty$ for some $\psi$) and this does not depend on the chart at infinity, so we actually obtain an element of $ (\nn_{\rl H^n})^*/O^+(n,1)  \cup \set{\infty}$.

\begin{rem}
The standard positive mass theorem, about asymptotically Euclidean metrics, can be thought of in a similar way. Let $M$ be a spin asymptotically Euclidean manifold with non-negative scalar curvature. Witten's trick consists in trying to build parallel spinors $\tilde{\psi}$ on $M$ that are asymptotic to the constant spinors $\psi$ of the Euclidean space $\rl^n$. Now every constant one-form on $\rl^n$ can be written as $X \mapsto i (X \cdot \psi,\psi)$. So the rigidity part of the Euclidean positive mass theorem can be explained as follows : starting from a constant one-form $\alpha$ on $\rl^n$, we pick a constant spinor $\psi$ such that $\alpha(X) = (X \cdot \psi,\psi)$ ; an analytical argument (based on $\mu_g=0$) provides a parallel spinor $\tilde{\psi}$ on $M$ asymptotic to $\psi$, hence a parallel one-form $\tilde{\alpha}$, given by $\tilde{\alpha}(X) = (X \cdot \tilde{\psi},\tilde{\psi})$, that is asymptotic to $\alpha$ ; this yields a parallel trivialization of the cotangent bundle of $M$, so $M$ is flat and is therefore $\rl^n$, owing to its asymptotic shape. The mass is the obstruction to do this. It is a single number $\mu$, but if we wish to make it fit into our picture, we might as well interpret it as a linear functional on the space of parallel sections of the flat bundle $\ee_{\rl^n} := T^* \rl^n \oplus \rl$ : $(\alpha,u) \mapsto \mu u$. The bundle $\ee_M := T^*M \oplus \rl$, endowed with the Levi-Civita connection on $T^*M$ and the flat connection on the $\rl$-part, is of course flat if and only if $M$ is flat, so the formalism described above still works.
\end{rem}

\section{A complex hyperbolic connection.}

\subsection{The connection.}

Let $(M^m,g,J)$ be a K\"ahler manifold of \emph{complex} dimension $m$. We wish to introduce a ``complex hyperbolic connection''
characterizing the complex hyperbolic geometry, in complete analogy with the real hyperbolic connection $\nabla^{RH}$ described in section 1. This is given by the following definition. The K\"ahler form is denoted by $\Omega := g(.,J.)$. We will often identify vectors and covectors thanks to the metric $g$. With our convention, if $(e_1,J e_1,\dots,e_m,J e_m)$ is an orthonormal basis, then $\Omega =\subup{k=1}{m}{\sum}  Je_k \wedge e_k$.

\begin{defn}\label{connexionCH}
Let $\ee := \Lambda^2_J M \oplus T^*M \oplus \rl$ be the vector bundle obtained as the direct sum of the bundle $\Lambda^2_J M$ of $J$-invariant $2$-forms, of the cotangent bundle and of the trivial (real) line bundle. We endow $\ee$ with the connection $\nabla^{CH}$ defined by
$$
\nabla^{CH}_X 
\begin{pmatrix}
\xi \\ 
\alpha 
\\ u 
\end{pmatrix}
:= 
\begin{pmatrix}
\nabla^g_X \xi + \frac{1}{2} \left( X \wedge \alpha + J X \wedge J \alpha \right) \\
\nabla^g_X \alpha + 2\iota_X ( \xi + u \Omega ) \\
d_X u + J \alpha (X)
\end{pmatrix}
.$$
The connection $\nabla^{CH}$ preserves a pseudo-Riemannian structure $h$ on $\ee$, with signature $(m^2 +1, 2m)$ and given by
$
h(\xi,\alpha,u) = \abs{\xi}^2_g + u^2 - \frac{\abs{\alpha}^2_g}{2}.
$
\end{defn}

This connection $\nabla^{CH}$ is very related to the study of Hamiltonian two-forms in \cite{ACG}, where a similar but more sophisticated connection is introduced. Let us explain why this connection is natural. In analogy with section 1, we expect its curvature to measure the deviation from the complex hyperbolic geometry (in the spirit of Cartan's connections). In order to write down an explicit formula for the curvature, we introduce an algebraic operation : if  $X$ and $Y$ are two vectors and $\gamma$ is an exterior form, we set:
$$
C_{X,Y}(\gamma) := (X \wedge \iota_Y \gamma - Y \wedge \iota_X \gamma)
+ (JX \wedge \iota_{JY} \gamma - JY \wedge \iota_{JX} \gamma).
$$

\begin{prop}
The curvature of $\nabla^{CH}$ is given by 
$$
\Rm^{\nabla^{CH}}_{X,Y} 
=
\begin{pmatrix}
\Rm^{\nabla}_{X,Y} -  C_{X,Y}  & 0 & 0 \\ 
0 & \Rm^{\nabla}_{X,Y} - \left[ 2 \Omega(X,Y)J + C_{X,Y} \right] &0\\
0 &0 &0 \\  
\end{pmatrix}
.$$ 
It follows that $\Rm^{\nabla^{CH}}$ vanishes if and only if the holomorphic sectional curvature is $-4$. In other words, $(\ee M,\nabla^{CH})$ is flat if and only if the universal cover of $(M^m,g,J)$ is the complex hyperbolic space $\cx H^m$ of holomorphic sectional curvature $-4$. 
\end{prop}

The proof of this formula is a straightforward computation, which we omit. The link with the curvature of the complex hyperbolic space is explained in paragraph IX.7 of \cite{KN} (where the sign convention is the opposite of ours). 

\begin{rem}\label{glcxconn}
By changing $\nabla^g_X \alpha + 2\iota_X ( \xi + u \Omega )$ into $\nabla^g_X \alpha - 2 c \iota_X ( \xi + u \Omega )$ in the formula for the connection, it is possible to obtain a family of connections characterizing every constant holomorphic curvature $4c$ (in particular : complex projective spaces). 
\end{rem}

The parallel sections $(\xi,\alpha,u)$ for $\nabla^{CH}$ obey $\alpha = Jdu$ and $\xi = -\frac12 \nabla^g \alpha - u \Omega$
so that they are determined by their third component $u$, which satisfies the third order equation 
$$
\Forall{X}{TM} \quad \nabla^g_X \hess_g u = 2 du(X) + X \odot du + JX \odot Jdu, 
$$
where $a \odot b$ means $a \otimes b + b \otimes a$. We will denote by $\nn$ the space of parallel sections of $(\ee, \nabla^{CH})$.
We will also need the subspace 
\begin{equation}\label{nzero}
\nn_0 := \set{(\xi,\alpha,u) \in \nn \; / \; g(\xi,\Omega) = u  }
\end{equation}
whose relevance will be clear from paragraph \ref{modelsec}.

\begin{rem}
The complex projective analogue ($c = +1$ in remark \ref{glcxconn}) of this third order equation appears in Obata's work \cite{Oba} :
on a simply-connected manifold, it possesses a solution if and only if the manifold is the complex projective space of holomorphic sectional curvature $4$. 
\end{rem}

\subsection{The model case.}\label{modelsec}

We wish to describe the case where $M = \cx H^m$. We start with the complex vector space $\cx^{m+1}$, $m \geq 2$, endowed with the Hermitian form $<,>$ defined by
$$
<z,z> = \sum_{k=1}^{m} \abs{z_k}^2 - \abs{z_{m+1}}^2.
$$
We denote this space by $\cx^{m,1}$ or $\rl^{2m,2}$ (whose metric structure is preserved by the standard complex structure $J$). The level set $<z,z> = -1$, endowed with the restriction of $<,>$ is by definition the Anti-de-Sitter space, $AdS^{2m,1}$, a Lorentz manifold with constant sectional curvature $-1$, invariant under
the natural (scalar) action of $\sph^1$ on $\cx^{m+1}$. The complex hyperbolic space $\cx H^m$ 
is the quotient $AdS^{2m,1} / \sph^1$, endowed with the induced metric :
\[
\begin{array}{cccc}
AdS^{2m,1} &     \subset & \cx^{m,1} \\
\pi\downarrow &             &  \\
\cx H^m    &             &  \\
\end{array}
\]
Let $\nu$ be the position vector field in $\cx^{m+1}$, identified also with the dual one-form $<\nu,.>$. 
By definition of $AdS^{2m,1}$, $<\nu,\nu>=-1$ along $AdS^{2m,1}$ and the tangent bundle $TAdS^{2m,1}$ is exactly the orthogonal subspace to $\nu$ for $<,>$. Besides, the vector field $J \nu$ is tangent to the action of $\sph^1$ and obeys $<J\nu,J\nu>=-1$. If  $\pi$ is the projection of $AdS^{2m,1}$ onto $\cx H^m$, it follows that at each point of $AdS^{2m,1}$, $d\pi$ is an isometry between $\set{\nu,J\nu}^{\bot}$ and the tangent space of $\cx H^m$. 

To understand $\ee = \ee(\cx H^m)$, we pick a point $z$ in $AdS^{2m,1}$ and look at 
the map $\appli{\theta_z}{\ee_{\pi(z)} }{\left( \Lambda^2_{J} \rl^{2m,2} \right)_z}$ given by 
$$
\theta_z(\xi,\alpha,u) = (d\pi_z)^* \xi + u(\pi(z)) \; J\nu \wedge \nu + \frac{ \nu \wedge (d\pi_z)^* \alpha  + J\nu \wedge J(d\pi_z)^*\alpha }{2}.
$$
It is a $\rl$-linear isomorphism and we have 
$$
< \theta_z(\xi,\alpha,u), \theta_z(\xi,\alpha,u) >
= \abs{\xi}^2_{\pi(z)} + u(\pi(z))^2 - \frac{1}{2} \abs{\alpha}^2_{\pi(z)} = h_{_{\pi(z)}}(\xi,\alpha,\nu).
$$
Setting $\Psi_z := \theta_z^{-1}$, we obtain a bundle map $\Psi$ that is $\sph^1$-invariant and satisfies 
the commutative diagram 
\diagcommut{\Lambda^2_{J} \rl^{2m,2}}{\Psi}{\ee}{}{}{AdS^{2m,1}}{\pi}{\cx H^m}
In this way, in the case of $\cx H^m$, the metric $h$ preserved by $\nabla^{CH}$ simply comes from $<,>$, via $\Psi$. At the level of \emph{sections}, $\Psi$ yields an isomorphism 
$$
\Gamma(\ee ) \cong \Gamma\left( \Lambda^2_{J} \rl^{2m,2} \restric{AdS^{2m,1}} \right)^{\sph^1}
$$
where the right-hand side denotes the $\sph^1$-invariant sections of the bundle obtained by restricting $\Lambda^2_{J} \rl^{2m,2}$ to $AdS^{2m,1}$. This isomorphism identifies $\sigma = (\xi,\alpha,u) \in \Gamma(\ee)$ with the $\sph^1$-invariant $2$-form 
$$
\Psi^*\sigma = \pi^* \xi + \pi^*u \; J\nu \wedge \nu +  \frac{\nu \wedge \pi^* \alpha + J\nu \wedge J\pi^*\alpha}{2}.
$$
The bundle $\Lambda^2_{J} \rl^{2m,2} \restric{AdS^{2m,1}}$ carries a natural connection $D$, inherited from the standard flat connection on $\rl^{2m,2}$. The real and imaginary parts of the differential forms $dz_k\wedge d\bar{z}_l$, restricted to  $AdS^{2m,1}$, trivialize the bundle $\Lambda^2_{J} \rl^{2m,2} \restric{AdS^{2m,1}}$, are $\sph^1$-invariant and $D$-parallel. $D$ therefore induces, via $\Psi$, a flat connection on $\ee $ and the reader might expect the following result. 

\begin{prop}\label{connexions}
The morphism $\Psi$ identifies $D$ and $\nabla^{CH}$ : for every horizontal vector field $H$ on $AdS^{2m,1}$ and for every section $\sigma$ of $\ee $, $\nabla^{CH}_{\pi_*H} \sigma = \Psi_* \; D_{H} \; \Psi^*\sigma$. 
\end{prop}

In this statement, ``horizontal'' means ``orthogonal to both $\nu$ and $J\nu$''. The proof of this is a direct computation, 
involving only two obvious facts : $D\nu$ is the identity at each point and the Levi-Civita connection of the hyperbolic space is induced by $D$. 

Note it is very important to work here with $J$-invariant $2$-forms. For instance, there is no flat $\sph^1$-invariant trivialization of the bundle $\Lambda^1 \rl^{2m,2} \restric{AdS^{2m,1}}$. 

We therefore obtain an isomorphism $\nn_{\cx H^m} \cong \Lambda^2_{J} \rl^{2m,2}$. Let $\Lambda^2_{J,0} \rl^{2m,2}$
be the subspace of primitive $J$-invariant $2$-forms on $\rl^{2m,2}$. Since the K\"ahler form $\omega$ of $\cx^{m,1}$, 
restricted to $AdS^{2m,1}$, is given by $\omega = - J\nu \wedge \nu + \pi^*\Omega$, where $\Omega$ is the K\"ahler form of $\cx H^m$, we also have (cf. (\ref{nzero})):
$$
\nn_{\cx H^m,0} \cong \Lambda^2_{J,0} \rl^{2m,2}.
$$

\subsection{K\"ahlerian Killing spinors.}

Let $(M^m,g,J)$ be a spin K\"ahler manifold. As in the hyperbolic space, a class of special spinors somehow characterizes 
the complex hyperbolic space. It turns out that there is a dimension issue, so we first assume the complex dimension $m$ is odd : $m=2l-1$. The even dimensional case will be discussed in section 4.

The spinor bundle $\Sigma$ decomposes into the orthogonal sum of the eigenspaces of the natural action of the K\"ahler form $\Omega=g(.,J.)$ : $\displaystyle{\Sigma = \bigoplus_{k=0}^m \Sigma_k}$, where $\Sigma_k$ corresponds to the eigenvalue $i (m-2k)$
(\cite{Kir} for instance). We will write $\pi^\Omega_k (\psi)$ or simply $\psi_k$ for the $k$th component of the spinor $\psi$ in this decomposition. From another point of view, spinors may be seen as twisted forms :
$\displaystyle{
\Sigma = \bigoplus_{k=0}^m \Lambda^{0,k} \otimes \sqrt{\Lambda^{m,0}}
}$
and then ${\Sigma_k=\Lambda^{0,k} \otimes \sqrt{\Lambda^{m,0}}}$. Through this identification, Clifford product ``$\cdot$'' is merely $\sqrt{2}$ times the difference between exterior product and interior product ($(1,0)$-vectors are identified with $(0,1)$-covectors by the Hermitian inner product). As a consequence, 
\begin{equation}\label{degre}
T^{1,0} \cdot \Sigma_k \subset \Sigma_{k+1} \quad \text{and} \quad T^{0,1} \cdot \Sigma_k \subset \Sigma_{k-1}.
\end{equation}

The following formula defines a connection $\hat{\nabla}$ on the vector bundle $\Sigma$ : 
\begin{eqnarray*}
\hat{\nabla}_X \psi &=& \nabla_X \psi + i c(X^{1,0}) \psi_{l-1} + i c(X^{0,1}) \psi_{l} \\
&=& \nabla_X \psi + \dfrac{i}{2} c(X) (\psi_{l-1} + \psi_{l}) + \dfrac{1}{2} c(JX) (\psi_{l-1} - \psi_{l}).
\end{eqnarray*}
The sections of $\Sigma_{l-1} \oplus \Sigma_{l}$ that are parallel for the connection $\hat{\nabla}$ are called imaginary \emph{K\"ahlerian Killing spinors}. The dimension of the space $\kk$ of K\"ahlerian Killing spinors is therefore at most $C_{2l}^{l}$. This bound is attained on $\cx H^{m}$, as we will see in details after corollary \ref{Qmap}. First, we make a number of general useful remarks.  

For future reference, let us introduce the (modified) Dirac operator $\hat{\dirac}$ which is naturally associated with the connection $\hat{\nabla}$ : if $(\epsilon_j)_j$ is any $g$-orthonormal frame, it is given by  
$$
\hat{\dirac} \psi := \sum_{j=1}^{2m} c(\epsilon_j) \hat{\nabla}_{\epsilon_j} \psi.
$$ 
In \cite{Her}, another modified Dirac operator $\dd$ is used, for analytical reasons : 
$$
\dd \psi := \hat{\dirac} \psi - i (m+1) (\psi - \psi_{l-1} - \psi_{l})
= \dirac \psi - i (m+1) \psi,
$$
where $\dirac$ is the standard Dirac operator. It should be noticed that K\"ahlerian Killing spinors are canceled by both $\hat{\dirac}$ and $\dd$.

\vskip 0.5cm

Given a spinor $\phi = \phi_{l-1}+ \phi_{l} \in \Sigma_{l-1} \oplus \Sigma_{l}$, we will use the notation 
$\tilde{\phi} := \phi_{l-1} - \phi_{l}$ ; $\tilde{\phi}$ is $(-1)^{l+1}$ times the conjugate $\phi_+ - \phi_-$ with respect 
to the usual decomposition into half-spinors. The K\"ahlerian Killing equation can then be written 
$$
\nabla_X \phi = -\frac{i}{2} X \cdot \phi - \frac{1}{2} (JX) \cdot \tilde{\phi}
\quad \text{ or } \quad 
\nabla_X \tilde{\phi} = \frac{i}{2} X \cdot \tilde{\phi} + \frac{1}{2} (JX) \cdot \phi.
$$
The following lemma is easy but very useful. 
\begin{lem}\label{trucs}
If $\phi$ is a K\"ahlerian Killing spinor then
\begin{eqnarray*}
\left(X\cdot\phi,\phi\right) &=& - \left(X\cdot \tilde{\phi},\tilde{\phi} \right) =2i \pim \left(X^{1,0}\cdot \phi_{l-1},\phi_l\right) \\
(X\cdot\tilde{\phi},\phi) &=& - \left(X\cdot \phi,\tilde{\phi} \right) = 2 \pre \left(X^{1,0}\cdot \phi_{l-1},\phi_l\right).
\end{eqnarray*}
In particular, $( (JX) \cdot \phi, \phi) = i (X \cdot \tilde{\phi}, \phi)$.
\end{lem}

\proof
The first statement follows from the computation 
\begin{eqnarray*}
	\left(X\cdot\phi,\phi\right)&=& \left(X\cdot(\phi_l +\phi_{l-1}),\phi_l +\phi_{l-1} \right) 
	= \left(X\cdot \phi_l, \phi_{l-1}\right) + \left(X\cdot \phi_{l-1},\phi_l\right) \\
	&=& 2i \pim \left(X\cdot \phi_{l-1},\phi_l\right) 
	= 2i \pim \left(X^{1,0}\cdot \phi_{l-1},\phi_l\right)
\end{eqnarray*}
and the other ones are similar.
\endproof

In complete analogy with the real hyperbolic case, the squared norm $\abs{\phi}^2$ of a K\"ahlerian Killing spinor $\phi$ will
induce a $\nabla^{CH}$-parallel section of $\ee$. To see this, we first compute two derivatives of this function.

\begin{lem}\label{der}
Any K\"ahlerian Killing spinor $\phi$ obeys 
\begin{eqnarray}
d \abs{\phi}^2(X) &=& -2i (X \cdot \phi, \phi) \label{der1} \\
\nabla_X d \abs{\phi}^2 (Y) &=& 2(X,Y) \abs{\phi}^2 - 2 \pim \left( Y \cdot (JX) \cdot \tilde{\phi}, \phi \right) \label{der2}.
\end{eqnarray}
\end{lem}

\proof
The definition of K\"ahlerian Killing spinors readily yields
$$
\d \left|\phi \right|^2 (X)
= -\dfrac{i}{2} (X \cdot \phi, \phi) - \dfrac{1}{2} (JX \cdot \tilde{\phi}, \phi)
+\dfrac{i}{2} (\phi, X \cdot \phi) - \dfrac{1}{2} (\phi, JX \cdot \tilde{\phi}),
$$
which leads to $d \abs{\phi}^2(X) = -2i (X \cdot \phi, \phi)$, thanks to lemma \ref{trucs}. Differentiating once more,
we find
\begin{eqnarray*}
\nabla_X d \abs{\phi}^2 (Y) 
&=& -2i (Y \cdot \nabla_X \phi,\phi) -2i (Y \cdot \phi, \nabla_X \phi) \\
&=& - (Y \cdot X \cdot \phi,\phi) + i (Y \cdot JX \cdot \tilde{\phi} ,\phi) 
+ (Y \cdot \phi, X \cdot \phi) + i (Y \cdot \phi, JX \cdot \tilde{\phi} )
\end{eqnarray*}
and the result follows from elementary properties of the Clifford product. 
\endproof

At this point, it is useful to introduce the following notations : for every K\"ahlerian Killing spinor
$\phi$, we define : 
$$
u_\phi := \abs{\phi}^2, \quad \quad \alpha_\phi := Jdu_\phi, \quad \quad
\xi_\phi(X,Y) :=  \pim (X \cdot Y \cdot \tilde{\phi}, \phi). 
$$

\begin{lem}\label{algxi}
$\xi_\phi$ is a $J$-invariant two-form satisfying :
\begin{eqnarray*}
\nabla_Z \alpha_\phi  &=& - 2 \iota_Z ( \xi_\phi + u_\phi \Omega) , \\
(\xi_\phi,\Omega) &=& u_\phi, \\
\nabla_Z \xi_\phi &=& - \frac{1}{2} \left( Z \wedge \alpha_\phi + J Z \wedge J\alpha_\phi \right) .
\end{eqnarray*}
\end{lem}

\proof
First, it is easy to check that $\xi_\phi$ is skewsymmetric, as a consequence of the definition of $\tilde{\phi}$. 
The $J$-invariance of $\xi_\phi$ stems from the following reformulation of its definition :
$$
\xi_\phi(X,Y) :=  \pim (X^{1,0} \cdot Y^{0,1} \cdot \tilde{\phi}, \phi) +\pim (X^{0,1} \cdot Y^{1,0} \cdot \tilde{\phi}, \phi). 
$$
Then formula (\ref{der2}) readily yields  
$$
\nabla_X \alpha_\phi (Y) = -2  (X,JY) u_\phi + 2 \xi_\phi(JY,JX) = -2 u_\phi \; \Omega(X,Y) - 2 \xi_\phi(X,Y),
$$ 
which justifies the first formula.

Let $(e_1,Je_1,\dots,e_m,Je_m)$ be an orthonormal basis. Then we have
$$
(\xi_\phi,\Omega) = \sum_{k=1}^m \xi_{\phi}(J e_k, e_k) = \sum_{k=1}^m \pim (J e_k \cdot e_k \cdot \tilde{\phi}, \phi)
= \pim (\Omega \cdot \tilde{\phi}, \phi).
$$
Now we use the spectral decomposition of the action of the K\"ahler form :
$$
(\Omega \cdot \tilde{\phi},\phi) = (\Omega \cdot \phi_{l-1} - \Omega \cdot \phi_{l}, \phi) 
= ( i \phi_{l-1} + i \phi_{l}, \phi)  = i \abs{\phi}^2.
$$
This ensures $(\xi_\phi,\Omega) = \abs{\phi}^2 = u_\phi$. 

Finally, to obtain the third equation, 
we introduce $\theta(X,Y) :=  2  (X \cdot Y \cdot \tilde{\phi}, \phi)$ and differentiate : 
\begin{eqnarray*}
& & \nabla_Z \theta(X,Y) \\
&=&  2(X \cdot Y \cdot \nabla_Z\tilde{\phi}, \phi) + 2(X \cdot Y \cdot \tilde{\phi}, \nabla_Z\phi) \\
&=&  i (X \cdot Y \cdot Z \cdot \tilde{\phi}, \phi) + (X \cdot Y \cdot (JZ) \cdot \phi, \phi) 
 + i (X \cdot Y \cdot \tilde{\phi}, Z \cdot \phi) -  (X \cdot Y \cdot \tilde{\phi}, (JZ) \cdot \tilde{\phi}) \\
&=&  i ( [X \cdot Y \cdot Z - Z \cdot X \cdot Y] \cdot \tilde{\phi}, \phi) + ( [X \cdot Y \cdot (JZ) - (JZ) \cdot X \cdot Y] \cdot \phi, \phi) . 
\end{eqnarray*}
The identity $ABC -CAB = 2 (A,C) B - 2(B,C) A$ in Clifford algebra leads to
\begin{eqnarray*}
& & \nabla_Z \theta(X,Y) \\
&=&  2i (X,Z) (Y \cdot \tilde{\phi},\phi) -2i (Y,Z) ( X \cdot \tilde{\phi}, \phi) 
+ 2 (X,JZ) (Y \cdot \phi,\phi) -2 (Y,JZ) ( X \cdot \phi, \phi). 
\end{eqnarray*}
Lemmata \ref{trucs} and \ref{der} then yield : 
\begin{eqnarray*}
& & \nabla_Z \theta(X,Y) \\
&=&  i (X,Z) du_\phi(JY) -i (Y,Z) du_\phi(JX) + i (X,JZ) du_\phi(Y) - i (Y,JZ) du_\phi(X) \\ 
&=& -i (X,Z) \alpha_\phi(Y) + i (Y,Z) \alpha_\phi(X) - i (X,JZ) J\alpha_\phi(Y) + i (Y,JZ) J\alpha_\phi(X)
\end{eqnarray*}
and the result follows. 
\endproof

\vskip 0.5cm

The computations above result in the following proposition (recall(\ref{nzero})).

\begin{cor}\label{Qmap}
If $\phi$ is a K\"ahlerian Killing spinor, then $Q(\phi) := (\xi_\phi,\alpha_\phi,u_\phi)$ is an element of $\nn_0$.
So we have a map $\appli{Q}{\kk}{\nn_0}$.
\end{cor}

When $m=2l-1$ for some integer $l \geq 2$, the complex hyperbolic space $\cx H^m$ is known to carry a space of K\"ahlerian Killing spinors of maximal dimension, $C^l_{2l}$ \cite{Kir}. Where do they come from ? Briefly, constant spinors on $\cx^{m,1}$ admit a restriction as spinors on $AdS^{2m,1}$ (cf. Lemma 3 in \cite{Bau}). These are the so-called imaginary Killing spinors ; they trivialize the spinor bundle of $AdS^{2m,1}$. Among them, thanks to the parity of $m$, some are $\sph^1$-invariant (it can be seen on the graduation coming from the K\"ahler structure of  $\cx^{m,1}$) and admit a projection into spinors along $\cx H^m$ (\cite{Mor}). These ``projected'' spinors are exactly the K\"ahlerian Killing spinors.

To be more explicit, we can adapt the computations of \cite{Kir}. The manifold $\cx H^m$ carries global coordinates $w_1, \dots, w_m$, such that 
$\pi^* w_k = \frac{z_k}{z_{m+1}}$ (along $AdS^{2m,1}$, where $z_{m+1}$ does not vanish). These coordinates induce a trivialization of the canonical bundle, by $\sqrt{dw}:=\sqrt{dw_1 \wedge \dots \wedge dw_m}$. We will use multi-index $a = (a_1,\dots,a_k)$ with $1 \leq a_1 < \dots < a_k \leq m$ and  set $dz_a := dz_{a_1} \wedge \dots \wedge dz_{a_k}$. The computations of \cite{Kir} say that, if $a$ is a multi-index of length $l-1$, the spinor $\varphi^a = \varphi^a_{l-1} + \varphi^a_{l}$ defined by
$$
\varphi^a_{l-1} = c(l) \frac{d\bar{w}_{a}}{(1-|w|^2)^l} \otimes \sqrt{dw} 
\quad
\text{ and }
\quad
\varphi^a_{l} = \frac{c(l)}{2i l} \bar{\partial} \left( \frac{d\bar{w}_{a} }{(1-|w|^2)^l} \right) \otimes \sqrt{dw} 
$$
is a K\"ahlerian Killing spinor. In this expression, $c(l)$ is a normalization constant, which we choose to be $c(l)= \sqrt{2}^{\frac52 - 3l}$. Another family of K\"ahlerian Killing spinors is described
by spinors $\breve{\varphi}^b = \breve{\varphi}^b_{l-1} + \breve{\varphi}^b_{l}$ where $b$ is a multi-index of length $l$ and 
$$
\breve{\varphi}^b_{l-1} =c(l)  \frac{ \iota_{\bar{R}} d\bar{w}_{b} }{(1-|w|^2)^l} \otimes \sqrt{dw} 
\quad
\text{ and }
\quad
\breve{\varphi}^b_{l} = \frac{c(l)}{2i l} \bar{\partial} \left(  \frac{\iota_{\bar{R}} d\bar{w}_{b}}{(1-|w|^2)^l} \right) \otimes \sqrt{dw} 
$$
where $R=\sum_k w_k \partial_{w_k}$. These families together form a basis for the K\"ahlerian Killing spinors of $\cx H^m$. What we are interested in here is their squared norms. Let us introduce the notation $\breve{a}$ for the multi-index that is complementary to $a$ (namely, $a$ and $\breve{a}$ have no common index and the sum of their lengths is $m$), and also the notation $\abs{w_a}^2 := \abs{w_{a_1}}^2 + \dots + \abs{w_{a_{k}}}^2$. An adaptation of the computations at the end of \cite{Kir} yields
$$
\abs{\varphi^a_{l-1}}^2 = |\breve{\varphi}^{\breve{a}}_{l}|^2 = \frac{1 - \abs{w_a}^2}{1 - \abs{w}^2}
\quad
\text{ and }
\quad
\abs{\varphi^a_{l}}^2 = |\breve{\varphi}^{\breve{a}}_{l-1}|^2 = \frac{\abs{w_{\breve{a}}}^2}{1 - \abs{w}^2},
$$ 
hence 
$$
\abs{\varphi^a}^2 = |\breve{\varphi}^{\breve{a}}|^2 = \frac{1 - \abs{w_a}^2 + \abs{w_{\breve{a}}}^2}{1 - \abs{w}^2}.
$$
It follows that $\pi^* u_{\varphi^a} = - \abs{z_a}^2 + \abs{z_{\breve{a}}}^2 + |z_{m+1}|^2$ (we forget the $\breve{\varphi}$'s since they do not yield new squared norms). Now $\alpha_{\varphi^a}$ and $\xi_{\varphi^a}$ are determined by $u_{\varphi^a}$ (since these are the components of an element of $\nn$). Denoting the real part of $z_k$ by $x_k$, one can check that the corresponding element of $\nn_0 (\cx H^m) \cong \Lambda^2_{J,0} \rl^{2m,2}$ is
$$
\beta_{\varphi^a} = -\sum^{l-1}_{j=1} Jdx_{a_j} \wedge dx_{a_j} +\sum^{l}_{j=1} Jdx_{\breve{a}_j} \wedge dx_{\breve{a}_j} + Jdx_{m+1} \wedge dx_{m+1}
$$
(To see this, it is sufficient to compute the scalar product of the right-hand side with $J\nu \wedge \nu$ and observe that 
it coincides with $\pi^* u_{\varphi^a}$.) In particular, $<\beta_{\varphi^a},\beta_{\varphi^a}> = m+1$.  All these two-forms 
belong to the same orbit of the isometric and holomorphic action of $PU(m,1)$. We will denote this orbit by $\nn_{0,\cx H^m}^+$.
Note that, since $\Lambda^2_{J,0} \rl^{2m,2}$ is an irreducible representation of $PU(m,1)$, the linear span of $\nn_{0,\cx H^m}^+$
is the whole  $\Lambda^2_{J,0} \rl^{2m,2}$.

\section{Toward a mass.} 

\subsection{The ``mass integral'' at infinity.}

Let us give a precise definition for the class of manifolds we are interested in.

\begin{defn}\label{defasympcx}
A complete K\"ahler manifold $(M^m,g,J)$ is called \emph{asymptotically complex hyperbolic} if there is a compact subset $K$ of $M$ and a ball $B$ in $\cx H^m$ such that:
\begin{enumerate}
	\item[(i)]  $(M \backslash K,J)$ is biholomorphic to $\cx H^m \backslash B$ and,
 \item[(ii)]  through this identification, $e^{a r} (g - g_{\cx H^m})$ is bounded in $C^{1,\alpha}$, with respect to the complex hyperbolic metric $g_{\cx H^m}$. Here, $r$ denotes the distance to some point in $\cx H^m$ and we assume $a > m + \frac12$.
\end{enumerate}
\end{defn}

In \cite{Her}, it is only assumed that the complex structure $J$ of $M$ is asymptotic to the complex structure $J_0$ of $\cx H^m$ (instead of $J=J_0$, as in our definition). It turns out that, if $J$ is asymptotic to $J_0$, then they are related by a biholomorphism, defined outside a compact set and asymptotic to the identity. To justify this, let us see $J$ and $J_0$ as two complex structures on a neighborhood $U$ of $\sph^{2m-1}$ in the unit ball of $\cx^m$. Under our assumptions, they induce the \emph{same} CR structure on the unit sphere, the standard CR structure of $\sph^{2m-1}$. The restriction of any $J_0$-holomorphic coordinate $z_k$ to the sphere is a $CR$ function. Since the standard sphere is strictly-pseudo-convex, Lewy's extension theorem makes it possible to extend this function into a $J$-holomorphic function $w_k$ on $U$ (shrinking $U$ if necessary). This yields a holomorphic map $w=(w_1,\dots,w_m)$ from $\Omega$ to $\cx^m$. Since $z=(z_1,\dots,z_m)$ is a diffeomorphism onto its range and coincide with $w$ on the sphere, we may shrink $U$ to ensure $w$ is also a diffeomorphism onto its range, hence a biholomorphism between two neighborhoods of $\sph^{2m-1}$ in the unit ball of $\cx^m$, one endowed with $J$ and the other one with the standard complex structure. The promised biholomorphism is $z^{-1} \circ w$.

In this section, we assume $(M,g,J)$ is an asymptotically complex hyperbolic manifold of odd complex dimension $m=2l-1$ and with 
scalar curvature bounded from below by $-4m(m+1)$ (the scalar curvature of $\cx H^m$). In this setting, the spinor bundle $\Sigma\restric{M \backslash K}$ can be identified with the spinor bundle of $\cx H^{m} \backslash B$. We may therefore extend any K\"ahlerian Killing spinor $\phi$ on $\cx H^{m} \backslash B$ into a spinor $\phi$ on $M$. Our aim is to understand to what extent we can make it into a K\"ahlerian Killing spinor. In the spirit of \cite{Wit}, we first check that we can choose an extension in the kernel of a Dirac operator. 

\begin{lem}\label{analyse}
There is a smooth spinor $\psi := \phi + \phi_{err}$ such that $\dd \psi= 0$ and $\phi_{err}$ decays as $e^{-b r}$ in $C^{1,\alpha}$
with $b > m+\frac12$.
\end{lem}

\proof
We first observe that 
$$
\dd \dd^* = \dirac^2 + (m+1)^2 = \nabla^*\nabla + \frac14 \scal + (m+1)^2 \geq \nabla^*\nabla + m+1.
$$
A slight adaptation of Proposition I.3.5 in \cite{Biq} then shows that the operator
$$
\appli{\dd \dd^*}{e^{b r} C^{2,\alpha}}{e^{b r} C^{0,\alpha}}
$$ 
is an isomorphism for every $b$ such that
$
m - \sqrt{m^2 + m+1} < b < m + \sqrt{m^2 + m+1}.
$
Since $\dd \phi = (\dd_g - \dd_{\cx H^m}) \phi = \oo(e^{-(a-1)r})$, with $a > m + \frac12$, we can pick a $b > m+\frac12$
such that the equation $\dd \dd^* \sigma = - \dd \phi_0$ admits a solution $\sigma$ in $e^{b r} C^{2,\alpha}$. Then  
$\phi_{err} := \dd^* \sigma$ is convenient.
\endproof

We then invoke a Weitzenb\"ock formula, proved in paragraph 3 of \cite{Her} (modulo two misprints, indeed: in the formula stated, an $i$ should be added at the second and third lines and the coefficient $m-2q$ and the fourth one should be replaced by $2(m-q)$) :
\begin{eqnarray*}
\int_{S_R} * \zeta_{\psi,\psi} 
=& &\int_{B_R} \abs{\hat{\nabla}\psi}^2 + \dfrac{1}{4} \int_{B_R} \left(\scal + 4m(m+1) \right) \abs{\psi}^2 \\
 &+& (m+1) \int_{B_R} \left(\abs{\psi}^2 - \abs{\pi^\Omega_{l-1}\psi}^2 - \abs{\pi_{l}^\Omega\psi}^2 \right),
\end{eqnarray*}
where $S_R$ denotes the sphere $\set{r=R}$, bounding the domain $B_R$, and $\zeta_{\sigma,\tau}$ is the $1$-form defined by
$$
\zeta_{\sigma,\tau}(X) = (\hat{\nabla}_X \sigma + c(X) \dd\sigma, \tau ).
$$
In view of this formula, the obstruction for $\psi$ to be a K\"ahlerian Killing spinor is precisely the ``mass integral at infinity'' $\lim_{R \to \infty} \int_{S_R} * \zeta_{\psi,\psi}$, which is a well defined element of $[0,+\infty]$, because the integrand on the right-hand side is non-negative.

\begin{lem}
$\displaystyle{
\lim_{R \to \infty} \int_{S_R} * \zeta_{\psi,\psi} = \lim_{R \to \infty} \int_{S_R} * \zeta_{\phi,\phi}.
}$ 
\end{lem}

\proof
Since $\pre \zeta_{\sigma,\tau}$ is symmetric up to a divergence term (as noticed in \cite{Her}, p.651),
we only need to check that 
$$
\lim_{R \to \infty} \int_{S_R} * (\zeta_{\phi,\phi_{err}} + \zeta_{\phi_{err},\phi_{err}} ) =0.
$$
This follows from the following estimates : $\vol S_R = \oo(e^{2m r})$,  
$\hat{\nabla} \phi = \oo(e^{(1-a) r})$ (beware $\phi$ grows in $e^{r}$), $\phi_{err} = \oo(e^{-b r})$, 
$\hat{\nabla} \phi_{err} = \oo(e^{-b r})$ 
(cf. lemma \ref{analyse}), with $a > m + \frac12$ and $b > m + \frac12$.
\endproof

In order to compare the metrics $g$ et $g_0 := g_{\cx H^m}$, we introduce the symmetric endomorphism $A$ such that $g_0 = g(A., A.)$. Since $A$ maps $g_0$-orthonormal frames to $g$-orthonormal frames, it identifies the spinor bundles defined with $g$ and $g_0$ (cf. \cite{CH} for instance). The associated Clifford products $c_g$ and $c_{g_0}$ are related by the formula ${c_g(AX)\sigma = c_{g_0}(X)\sigma}$. Note we will 
also write $X \cdot$ for $c_{g_0}(X)$ (and not for $c_g(X)$). For the sake of efficiency, we will write $u \approx v$ when $u-v = o(e^{-2mr})$ ; the terms we neglect in this way will indeed not contribute to the integral at infinity. Before computing the ``mass integral at infinity'', we point out a few elementary facts.

\begin{lem}\label{commut}
$\displaystyle{ A^{-1} J A \approx J}$.
\end{lem}

\proof
The definition $g_0=g(A,A)$ and the compatibility of $J$ with $g_0$ yield the equality ${g(AJ,AJ)=g(A,A)}$. Since $A$ and $J$ are respectively $g$-symmetric and $g$-antisymmetric, we deduce : $JA^2J = -A^2$. If $A=1+H$, this implies $ JHJ \approx -H$. Since $J^2=-1$, we obtain $ JH \approx H J$ and then $ J A \approx A J$, hence the result. 
\endproof

\begin{cor}\label{kahlcliff}
We have $\displaystyle{ c_g(\Omega) \approx c_{g_0}(\Omega_0)}$ and $\displaystyle{ \pi_k^{\Omega} \approx \pi_k^{\Omega_0}}$. 
\end{cor}

\proof
Given a $g_0$-orthonormal basis $(e_1,J e_1,\dots,e_{m}, J e_{m})$,  the Clifford action of the K\"ahler form $\Omega_0$ reads 
$c_{g_0}(\Omega_0) = \sum_{k} J e_k \cdot e_k \cdot$ while the K\"ahler form $\Omega$ of $g$ acts by 
$$
c_g(\Omega) = \sum^{m}_{k=1} c_g (J A e_k) c_g(A e_k) = \sum^{m}_{k=1} (A^{-1}J Ae_k) \cdot e_k \cdot .
$$ 
The first statement is therefore a straightforward consequence
of lemma \ref{commut}. The second one follows from general considerations. We observe the skew-Hermitian endomorphisms $P := c_g(\Omega)$ and $P_0 := c_{g_0}(\Omega_0)$ act on each fiber of the spinor bundle with the same spectrum. If $\lambda$ is one of the eigenvalues, the corresponding spectral projectors $\Pi$ and $\Pi_0$ (for $P$ and $P_0$) obey the formulas 
$$
\Pi = \frac{1}{2\pi} \int_{C} (z-P)^{-1} dz \quad \text{and} \quad \Pi_0 = \frac{1}{2\pi} \int_{C} (z-P_0)^{-1} dz
$$ 
where $C$ is a circle in the complex plane, centered in $\lambda$ and with small radius $\delta$. We deduce
$$
\Pi -\Pi_0 = \frac{1}{2\pi} \int_{C} (z-P)^{-1} (P-P_0) (z-P_0)^{-1} dz
$$ 
and then 
$$
\abs{\Pi -\Pi_0} \leq \delta \delta^{-1} \abs{P-P_0} \delta^{-1} = \delta^{-1} \abs{P-P_0}.
$$ 
The result follows at once.
\endproof

The rest of this section is devoted to the proof of the following statement.

\begin{prop}\label{formmass}
The ``mass integral at infinity'' is
$$
\lim_{R \to \infty} \int_{S_R} * \zeta_{\psi,\psi} = \lim_{R \to \infty} \int_{S_R} * \left( -\dfrac{1}{4} \left( d\tr_{g_0} g + \Div_{g_0} g \right) \abs{\phi}^2 
+\dfrac{1}{8} \tr_{g_0} (g-g_0) \; \d \left|\phi \right|^2 \right).
$$
\end{prop}

\proof
To begin with, in view of corollary \ref{kahlcliff}, we may write
\begin{eqnarray*}
\zeta_{\phi,\phi}(Y) 
&=& (\hat{\nabla}^g_Y \phi + c_g(Y) \hat{\dirac}\phi,\phi) 
- i (m+1) \left( c_g(Y) (1 - \pi^\Omega_{l-1} - \pi^\Omega_{l}) \phi , \phi \right) \\
&\approx& (\hat{\nabla}^g_Y \phi + c_g(Y) \hat{\dirac}\phi,\phi).
\end{eqnarray*}
Given a $g_0$-orthonormal frame $(e_1,\dots,e_{2m})$ and a $g_0$-unit vector $X$, since $\phi$ is a K\"ahlerian Killing 
spinor with respect to $(g_0,J_0)$, we may therefore write outside $K$ (as in \cite{CH,Min} for instance) :
\begin{eqnarray*}
\zeta_{\phi,\phi}(AX) 
&\approx& \dfrac{1}{2}\sum_{j=1}^{2m} ([c_g(AX), c_g(A e_j)] \hat{\nabla}^g_{A e_j} \phi,\phi) \\
&\approx& \dfrac{1}{2}\sum_{j=1}^{2m} ([c_g(AX), c_g(A e_j)] (\hat{\nabla}^g_{A e_j}-\hat{\nabla}^{g_0}_{A e_j}) \phi,\phi) \\
&\approx& \dfrac{1}{2}\sum_{j=1}^{2m} ([X \cdot, e_j \cdot] (\hat{\nabla}^g_{A e_j}-\hat{\nabla}^{g_0}_{A e_j}) \phi,\phi), 
\end{eqnarray*}
which expands into
\begin{eqnarray*}
& &\dfrac{1}{2}\sum_{j=1}^{2m} \left( [X \cdot, e_j \cdot] (\nabla^g_{A e_j}-\nabla^{g_0}_{A e_j}) \phi,\phi \right) \\
+ &\dfrac{i}{4}& \sum_{j=1}^{2m} \left( [X \cdot, e_j \cdot]  \Big{(} c_g(A e_j-iJA e_j) \pi_{l-1}^\Omega -   c_{g_0}(A e_j-iJA e_j) \pi_{l-1}^{\Omega_0} \Big{)} \phi ,\phi \right)\\
+ &\dfrac{i}{4}& \sum_{j=1}^{2m} \left( ( [X \cdot, e_j \cdot] \Big{(}  c_g(A e_j+iJA e_j) \pi_{l}^\Omega -  c_{g_0}(A e_j+iJA e_j) \pi_{l}^{\Omega_0}\Big{)} \phi ,\phi\right),
\end{eqnarray*}
that is
\begin{eqnarray*}
& &\dfrac{1}{2}\sum_{j=1}^{2m} \left( [X \cdot, e_j \cdot] (\nabla^g_{A e_j}-\nabla^{g_0}_{A e_j}) \phi,\phi \right) \\
+ &\dfrac{i}{4}& \sum_{j=1}^{2m} \left( [X \cdot, e_j \cdot]  \Big{(} ( e_j-iA^{-1}JA e_j)\cdot \pi_{l-1}^\Omega -   (A e_j-iJ A e_j) \cdot \pi_{l-1}^{\Omega_0} \Big{)} \phi ,\phi \right)\\
+ &\dfrac{i}{4}& \sum_{j=1}^{2m} \left( ( [X \cdot, e_j \cdot] \Big{(}  (e_j+iA^{-1}JA e_j) \cdot \pi_{l}^\Omega -  (A e_j+iJ A e_j) \cdot \pi_{l}^{\Omega_0}\Big{)} \phi ,\phi\right).
\end{eqnarray*}
In view of lemma \ref{commut} and corollary \ref{kahlcliff}, we are left with :
\begin{eqnarray*}
\zeta_{\phi,\phi}(AX) \approx
& &\dfrac{1}{2}\sum_{j=1}^{2m} \left( [X \cdot, e_j \cdot] (\nabla^g_{A e_j}-\nabla^{g_0}_{A e_j}) \phi,\phi \right) \\
+ &\dfrac{i}{4}& \sum_{j=1}^{2m} \left( [X \cdot, e_j \cdot]  ( (e_j-iJ e_j) -  (A e_j-iJA e_j) ) \cdot \pi_{l-1}^{\Omega_0}  \phi ,\phi \right)\\
+ &\dfrac{i}{4}& \sum_{j=1}^{2m} \left( [X \cdot, e_j \cdot] ( (e_j+iA e_j) - (A e_j+iJA e_j) ) \cdot \pi_{l}^{\Omega_0} \phi ,\phi\right).
\end{eqnarray*}
With $A=1+H$ we can there therefore write  $\zeta_{\phi,\phi}(AX) \approx I + II + III$, with : 
\begin{eqnarray*}
I &:=& \dfrac{1}{2}\sum_{j=1}^{2m} \left( [X \cdot, e_j \cdot] (\nabla^g_{A e_j}-\nabla^{g_0}_{A e_j}) \phi,\phi \right) \\
II &:=& - \dfrac{i}{4} \sum_{j=1}^{2m} \left( [X \cdot, e_j \cdot]  H e_j \cdot \phi ,\phi \right) \\
III &:=& - \dfrac{1}{4} \sum_{j=1}^{2m} \left( [X \cdot, e_j \cdot] JH e_j \cdot \tilde{\phi} ,\phi\right).
\end{eqnarray*}
The computation of the real part of the first term is classical (cf. \cite{CH} or lemma 10 in \cite{Min}, for instance) :
$$
\pre I \approx -\dfrac{1}{4} \left( d\tr_{g_0} g + \Div_{g_0} g \right)\abs{\phi}^2. 
$$
The second term is basically computed in \cite{CH}. Indeed, since $H$ is symmetric, we use the identity 
$[X\cdot ,e_j\cdot ]=2X_j + 2 X \cdot e_j \cdot$ to obtain
\begin{eqnarray*}
\pre II &=& -\dfrac{i}{2} \left(H X \cdot \phi ,\phi \right) 
- \dfrac{i}{2} \sum_{j,l=1}^{2m} H_{jl} \left(X \cdot e_j \cdot e_l \cdot \phi ,\phi \right)\\
&=& -\dfrac{i}{2} \left(H X \cdot \phi ,\phi \right) + \dfrac{i}{2} \tr H \left( X \cdot \phi ,\phi \right).
\end{eqnarray*}
In the same way, the third term can be written
\begin{eqnarray*}
III = - \dfrac{1}{2} \left( ( JH X \cdot \tilde{\phi} ,\phi\right)
- \dfrac{1}{2} \sum_{j=1}^{2m} \left(  X \cdot e_j \cdot JH e_j \cdot \tilde{\phi} ,\phi\right).
\end{eqnarray*}

\begin{lem}
$\displaystyle{\pre \sum_{j=1}^{2m} \left(  X \cdot e_j \cdot JH e_j \cdot \tilde{\phi} ,\phi\right) 
\approx -2 \left(  JH X \cdot \tilde{\phi} ,\phi\right). }$
\end{lem}

\proof
Let us set $M := JH$ and $M_{i j} := (e_i,M e_j)$, so that
$$
\sum_{j=1}^{2m} \left(  X \cdot e_j \cdot JH e_j \cdot \tilde{\phi} ,\phi\right)
= \sum_{j,k,p} M_{k j} X_p \left(  e_p \cdot e_j \cdot e_k \cdot \tilde{\phi} ,\phi\right).
$$
Lemma \ref{commut} ensures $J H \approx HJ$. Since $H$ is symmetric and $J$ antisymmetric, we deduce 
that $M$ is antisymmetric up to a negligible term. In particular, $M_{k j} \approx 0$ when $k=j$, hence 
$$
\sum_{j=1}^{2m} \left(  X \cdot e_j \cdot JH e_j \cdot \tilde{\phi} ,\phi\right)
\approx \sum_{j \not= k} \sum_p M_{k j} X_p \left(  e_p \cdot e_j \cdot e_k \cdot \tilde{\phi} ,\phi\right).
$$
Given three \emph{distinct} indices $j,k,p$ we consider the expression 
$$\left(  e_p \cdot e_j \cdot e_k \cdot \tilde{\phi} ,\phi\right)
= \left(  e_p \cdot e_j \cdot e_k \cdot (\phi_{l-1} - \phi_l) , (\phi_{l-1} + \phi_l) \right).
$$
Property (\ref{degre}) reduces it into
$$
\left(  e_p \cdot e_j \cdot e_k \cdot \tilde{\phi} ,\phi\right)
= \left(  e_p \cdot e_j \cdot e_k \cdot \phi_{l-1} , \phi_l \right)
- \left(  e_p \cdot e_j \cdot e_k \cdot \phi_l , \phi_{l-1} \right).
$$
and since the indices are distinct, this is imaginary. So 
\begin{eqnarray*}
& &\pre \sum_{j=1}^{2m} \left(  X \cdot e_j \cdot JH e_j \cdot \tilde{\phi} ,\phi\right) \\
&\approx& \pre \sum_{j \not= k} M_{k j} X_j \left(  e_j \cdot e_j \cdot e_k \cdot \tilde{\phi} ,\phi\right) 
+ \pre \sum_{j \not= k} M_{k j} X_k \left(  e_k \cdot e_j \cdot e_k \cdot \tilde{\phi} ,\phi\right) \\
&\approx& -\pre \sum_{j \not= k} M_{k j} X_j \left(  e_k \cdot \tilde{\phi} ,\phi\right) 
+ \pre \sum_{j \not= k} M_{k j} X_k \left(  e_j \cdot  \tilde{\phi} ,\phi\right) \\
&\approx& -2 \pre \sum_{j \not= k} M_{k j} X_j \left(  e_k \cdot \tilde{\phi} ,\phi\right) \\
&\approx& -2 \pre \left(  MX \cdot \tilde{\phi} ,\phi\right) \\
&\approx& -2  \left(  MX \cdot \tilde{\phi} ,\phi\right).
\end{eqnarray*}
\endproof
This lemma leads to
$$
\pre III \approx - \dfrac{1}{2} \left(  JH X \cdot \tilde{\phi} ,\phi\right)
+ \left(  JH X \cdot \tilde{\phi} ,\phi\right) = \dfrac{1}{2} \left(  JH X \cdot \tilde{\phi} ,\phi\right).
$$
Eventually, summing $I$, $II$ and $III$, we get : 
\begin{eqnarray*}
\pre \zeta_{\phi,\phi}(AX) \approx 
&-&\dfrac{1}{4} \left( d\tr_{g_0} g + \Div_{g_0} g \right)\abs{\phi}^2 \\
&-&\dfrac{i}{2} \left(H X \cdot \phi ,\phi \right) + \dfrac{i}{2} \tr H \left( X \cdot \phi ,\phi \right) \\
&+& \dfrac{1}{2} \left(  JH X \cdot \tilde{\phi} ,\phi \right) .
\end{eqnarray*}
Lemmata \ref{trucs} and \ref{der} simplify this into 
\begin{eqnarray*}
\pre \zeta_{\phi,\phi}(AX) 
&\approx& 
-\dfrac{1}{4} \left( d\tr_{g_0} g + \Div_{g_0} g \right)\abs{\phi}^2 
+ \dfrac{i}{2} \tr H \left( X \cdot \phi ,\phi \right) \\
&\approx& -\dfrac{1}{4} \left( d\tr_{g_0} g + \Div_{g_0} g \right)\abs{\phi}^2 
- \dfrac{1}{4} \tr H \; \d \left|\phi \right|^2 (X).
\end{eqnarray*}
Since $g-g_0 \approx -2 g_0(H.,.)$, we have $\tr (g-g_0) \approx -2 \tr H$, hence :
$$
\pre \zeta_{\phi,\phi}(AX) \approx -\dfrac{1}{4} \left( d\tr_{g_0} g + \Div_{g_0} g \right)(X)\abs{\phi}^2 
+\dfrac{1}{8} \tr_{g_0} (g-g_0) \; \d \left|\phi \right|^2 (X),
$$
which yields the formula of proposition \ref{formmass}.
\endproof

\subsection{The mass linear functional.}

We consider the formula 
\begin{equation}\label{theformula}
\mu_g (\xi,\alpha,u) = 
-\dfrac{1}{4} \lim_{R \to \infty} \int_{S_R} * \left( \left( d\tr_{g_0} g + \Div_{g_0} g \right) u 
+\dfrac{1}{2} \tr_{g_0} (g-g_0) \; J\alpha \right).
\end{equation}
If $(\xi,\alpha,u)$ belongs to $Q_{\cx H^m}(\kk_{\cx H^m})$, the considerations of the previous paragraph
imply that $\mu_g (\xi,\alpha,u)$ is an element of $[0,+\infty]$. If it is infinite at some point of 
$Q_{\cx H^m}(\kk_{\cx H^m})$, we decide that the mass is infinite : $\mu_g = \infty$. Otherwise, this defines a linear 
functional $\mu_g$ on the linear span of $Q_{\cx H^m}(\kk_{\cx H^m})  \subset \nn_{\cx H^m,0}$, 
which is $\nn_{\cx H^m,0}$, since
it is $PU(m,1)$-invariant and the action of this group on $\Lambda^2_{J,0} \rl^{2m,2}$ is irreducible. 

The previous paragraph also indicates that $\mu_g$ takes non-negative values on (the convex cone generated by) 
$Q_{\cx H^m}(\kk_{\cx H^m})$ ; since $Q_{\cx H^m}(\kk_{\cx H^m})$ is $PU(m,1)$-invariant and contains one 
element of $\nn_{\cx H^m,0}^+$ (from the explicit computations in 2.3), it contains $\nn_{\cx H^m,0}^+$. 
So $\mu_g$ is non-negative on $\nn_{\cx H^m,0}^+$.

Assume $\mu_g$ vanishes on $\nn_{\cx H^m,0}$. Then Lemma \ref{analyse}, coupled to the Bochner formula, ensures that for every $\psi$ 
in $\kk_{\cx H^m}$, there is a K\"ahlerian Killing spinor $\phi$ on $(M,g)$, namely $\phi \in \kk_g$, that is asymptotic to $\psi$.
In particular, for any element $\beta = Q_{\cx H^m}(\psi)$ of $\nn_{\cx H^m,0}^+$, there is an element $\tilde{\beta} = Q_g(\phi)$ 
of $\nn_{g,0}$ that is asymptotic to $\beta$. Since the linear span of $\nn_{\cx H^m,0}^+$ is the whole $\nn_{\cx H^m,0}$, we deduce that $\nn_{g,0}$ has maximal dimension, which implies that 
$\nn_g$ has maximal dimension (add the K\"ahler form) : $(M,g)$ is (locally) complex hyperbolic. In view of its asymptotic, it is bound to be $\cx H^m$. 

The holomorphic chart at infinity in Definition \ref{defasympcx} is of course not unique. Any two relevant charts $\Psi_1$ 
and $\Psi_2$ differ by a biholomorphism $f := \Psi_2 \circ \Psi_1^{-1}$ of $\cx H^m$ such that $f^* g_0$ is asymptotic to 
$g_0$. The map $f$ therefore induces a CR-automorphism of the sphere at infinity, so that $f$ is asymptotic to an element 
of $PU(m,1)$. As in the remark after Definition \ref{defasympcx}, we conclude $f$ is an element of $PU(m,1)$. So the model at infinity is unique up to its natural automorphism group $PU(m,1)$. If $\phi_1$ is a K\"ahlerian Killing spinor in the chart $\Psi_1$, then $\phi_2 := f_* \phi_1$ is a K\"ahlerian Killing spinor in the chart $\Psi_2$ (by naturality). We claim that
\begin{lem}
$\displaystyle{
\mu_g^{\Psi_1}(\beta_{\phi_1}) = \mu_g^{\Psi_2}(\beta_{\phi_2}).
}$
\end{lem}

\proof
Recall that in the notations introduced in the previous paragraph,
$$
\mu_g^{\Psi_1}(\beta_{\phi_1}) = \lim_{R \to \infty} \int_{S_R^{\Psi_1}} *\zeta_{\phi_1,\phi_1} 
= \lim_{R \to \infty} \int_{S_R^{\Psi_1}} *\zeta_{\psi_1,\psi_1}.
$$
This quantity is equal to the well-defined integral 
\begin{eqnarray*}
& & \int_M \abs{\hat{\nabla}\psi_1}^2 + \dfrac{1}{4} \int_M \left(\scal + 4m(m+1) \right) \abs{\psi_1}^2 \\
& & + (m+1) \int_M \left(\abs{\psi_1}^2 - \abs{\pi^\Omega_{l-1}\psi_1}^2 - \abs{\pi_{l}^\Omega\psi_1}^2 \right),
\end{eqnarray*}
so that we can compute it with another family of spheres :
$$
\mu_g^{\Psi_1}(\beta_{\phi_1}) = \lim_{R \to \infty} \int_{f(S_R^{\Psi_1})} *\zeta_{\psi_1,\psi_1} 
= \lim_{R \to \infty} \int_{f(S_R^{\Psi_1})} *\zeta_{\phi_1,\phi_1}.
$$
Since $f^* g_0 = g_0$ and $f(S_R^{\Psi_1}) = S_R^{\Psi_2}$, we therefore obtain
\begin{eqnarray*}
\mu_g^{\Psi_1}(\beta_{\phi_1}) 
&=& \lim_{R \to \infty} \int_{f(S_R^{\Psi_1})} *_{g_0} \left( -\dfrac{1}{4} \left( d\tr_{g_0} g + \Div_{g_0} g \right) \abs{ \phi_1}^2_{g_0}  +\dfrac{1}{8} \tr_{g_0} (g-g_0) \; d\abs{\phi_1}^2_{g_0} \right) \\
&=& \lim_{R \to \infty} \int_{S_R^{\Psi_2}} *_{g_0} \left( -\dfrac{1}{4} \left( d\tr_{g_0} g + \Div_{g_0} g \right) \abs{f^*\phi_2}^2_{g_0}  +\dfrac{1}{8} \tr_{g_0} (g-g_0) \; d\abs{f^*\phi_2}^2_{g_0} \right).\\
&=& \lim_{R \to \infty} \int_{S_R^{\Psi_2}} *_{g_0} \left( -\dfrac{1}{4} \left( d\tr_{g_0} g + \Div_{g_0} g \right) \abs{\phi_2}^2_{g_0}  +\dfrac{1}{8} \tr_{g_0} (g-g_0) \; d\abs{\phi_2}^2_{g_0} \right) \\
&=& \mu_g^{\Psi_2}(\beta_{\phi_2}).
\end{eqnarray*}
\endproof
In other words, changing the chart at infinity by an automorphism $f \in PU(m,1)$ results in turning $\mu_g$
to $f^*\mu_g$. So  $\mu_g$ is well-defined up to the natural action of $PU(m,1)$. We have proved the following result. 

\begin{thm}
Let $(M^m,g,J)$ be a spin asymptotically complex hyperbolic K\"ahler manifold with odd complex dimension
$\scal_g \geq \scal_{\cx H^m}$. Then the linear functional $\mu_g$ on $\nn_{\cx H^m,0}$ is well-defined up to
the natural action of $PU(m,1)$ ; it is non-negative on $\nn^+_{\cx H^m,0}$ and vanishes if and only if 
$(M^m,g,J)$ is the complex hyperbolic space.
\end{thm}

\section{The case of even-dimensional manifolds.}

\subsection{The twisted K\"ahlerian Killing spinors.}

To extend the ideas above to the even-dimensional case, it seems that we need K\"ahlerian Killing spinors
on even-dimensional complex hyperbolic spaces. Unfortunately, such K\"ahlerian Killing spinors do not exist. 
To overcome this cruel reality, we follow \cite{BH} and turn to the spin$^c$ realm. We refer to \cite{LM}
for basic definitions about spin$^c$ structures.

We consider a K\"ahler manifold $(M,g,J)$ with even complex dimension $m=2l$. As in \cite{BH},
we further assume that the cohomology class of $\frac{\Omega}{i\pi}$ is integral, i.e. in the image of 
$H^2(M,\ir) \to H^2(M,\rl)$. This determines a complex line bundle $L$ endowed with a Hermitian metric and a unitary connection with curvature $F = -2i \Omega$ ; the Chern class of $L$ is $c_1(L) = \frac{i}{2\pi} [F] = \frac1\pi [\Omega]$. We also assume that 
$(M,L)$ defines a spin$^c$ structure, in that the bundle $TM \otimes L$ admits a spin structure. We then introduce 
the corresponding spinor bundle $\Sigma^c$, endowed with a Clifford action $c$ and a connection $\nabla$. The K\"ahler 
form $\Omega$ acts on this bundle, with the eigenvalues $i(m-2k)$,
$0 \leq k \leq m$. The eigenspaces yield subbundles $\Sigma^c_k$. We may therefore define a connection 
$\hat{\nabla}$ on $\Sigma^c$ by requiring  that
$$
\hat{\nabla}_X \psi := \nabla_X \psi + i c(X^{1,0}) \psi_{l-1} + i c(X^{0,1}) \psi_{l}, 
$$
where $\psi_k$ is the component of $\psi$ in $\Sigma^c_k$. The parallel sections of $\Sigma_{l-1}^c \oplus \Sigma^c_{l}$ for this connection will be called twisted K\"ahlerian Killing spinors and the corresponding subspace will be denoted by $\kk^c$. 

Every computation of section 2.3 (but one) can be carried out with twisted K\"ahlerian Killing spinors, leading to the same formulas. In particular, 
one can define $u_\psi$, $\alpha_\psi$, $\xi_\psi$ for any $\psi$ in $\kk^c$ as in section 2.3 and this yields a map $Q^c \, : \kk^c \to \nn$. The only difference is that we do not get elements of $\nn_0$ : $(\xi_\psi,\Omega)$ is no longer $u_\psi$, basically because the eigenvalues of the K\"ahler form are now \emph{even} (compare with the proof of Lemma \ref{algxi}).

Let us describe the model case, where $(M,g,J)$ is the complex hyperbolic space $\cx H^m$, $m=2l$, with holomorphic sectional curvature $-4$. With this normalization, we have $Ric_{\cx H^m} = - 2(m+1) g_{\cx H^m}$. Since the Ricci form $\rho = \ric(.,J.)$ is $i$ times the curvature of the canonical line bundle, this implies : 
$$
c_1(\Lambda^{m,0}) = \frac{i}{2\pi} [-i \rho] 
= \frac{i}{2\pi} [ 2i (m+1) \Omega] = -\frac{m+1}{\pi} [\Omega] = -(m+1) c_1(L).
$$
So $L \cong \left( \Lambda^{m,0} \right)^{-\frac{1}{m+1}}$ 
in this case. It follows that
$$
\Sigma^c_k \cong \Sigma_k \otimes L \cong \Lambda^{0,k} \otimes \left( \Lambda^{m,0} \right)^{\frac12 - \frac{1}{m+1}}
\cong \Lambda^{0,k} \otimes \left( \Lambda^{m,0} \right)^{\frac{l}{2l+1}}.
$$
It turns out that the sections of $\kk^c$ trivialize the bundle $\Sigma_{l-1}^c \oplus \Sigma^c_{l}$, as noticed in  \cite{BH}. Indeed, using the same notations as in section 2.3, we can define two families of twisted K\"ahlerian Killing spinors in the following way.  
First, if $a$ is a multi-index of length $l-1$, we set  $\varphi^a = \varphi^a_{l-1} + \varphi^a_{l}$ with
$$
\varphi^a_{l-1} = c(l) \frac{d\bar{w}_{a}}{(1-|w|^2)^l} \otimes dw^{\frac{l}{2l+1}} 
\quad
\text{ and }
\quad
\varphi^a_{l} = \frac{c(l)}{2i l} \bar{\partial} \left( \frac{d\bar{w}_{a} }{(1-|w|^2)^l} \right) \otimes dw^{\frac{l}{2l+1}}. 
$$
The normalization we choose is  $c(l)= 2^{1-l - \frac{l^2}{2l+1}}$. We also introduce $\breve{\varphi}^b = \breve{\varphi}^b_{l-1} + \breve{\varphi}^b_{l}$ where $b$ is a multi-index of length $l$ and 
$$
\breve{\varphi}^b_{l-1} =c(l)  \frac{ \iota_{\bar{R}} d\bar{w}_{b} }{(1-|w|^2)^l} \otimes dw^{\frac{l}{2l+1}} 
\quad
\text{ and }
\quad
\breve{\varphi}^b_{l} = \frac{c(l)}{2i l} \bar{\partial} \left(  \frac{\iota_{\bar{R}} d\bar{w}_{b}}{(1-|w|^2)^l} \right) \otimes dw^{\frac{l}{2l+1}} 
$$
where $R=\sum_k w_k \partial_{w_k}$ as before. These families together form a basis for $\kk^c$ on $\cx H^m$. 
One can again compute the squared norms of these spinors, like in \cite{Kir}. If $a$ is a multi-index of length $l-1$
and $b$ a multi-index of length $l$, we have
$$
\abs{\varphi^a_{l-1}}^2 = \frac{1 - \abs{w_a}^2}{1 - \abs{w}^2},
\quad
\abs{\varphi^a_{l}}^2  = \frac{\abs{w_{\breve{a}}}^2}{1 - \abs{w}^2},
\quad
\abs{\breve{\varphi}^{b}_{l-1}}^2 = \frac{\abs{w_b}^2}{1 - \abs{w}^2},
\quad
\abs{\breve{\varphi}^{b}_{l}}^2 = \frac{1-\abs{w_{\breve{b}}}^2}{1 - \abs{w}^2},
$$ 
hence 
$$
\abs{\varphi^a}^2  = \frac{1 - \abs{w_a}^2 + \abs{w_{\breve{a}}}^2}{1 - \abs{w}^2}, \quad
\abs{\breve{\varphi}^b}^2  = \frac{1 -\abs{w_{\breve{b}}}^2  + \abs{w_{b}  }^2    }{1 - \abs{w}^2}.
$$
Exactly as in the odd dimensional case, we see that the corresponding elements of $\nn (\cx H^m) \cong \Lambda^2_{J} \rl^{2m,2}$ are
\begin{eqnarray*}
Q^c(\varphi^a) = \beta_{\varphi^a} &=& -\sum^{l-1}_{j=1} Jdx_{a_j} \wedge dx_{a_j} +\sum^{l+1}_{j=1} Jdx_{\breve{a}_j} \wedge dx_{\breve{a}_j} + Jdx_{m+1} \wedge dx_{m+1}, \\
Q^c(\breve{\varphi}^b) = \beta_{\breve{\varphi}^b} &=& -\sum^{l}_{j=1} Jdx_{b_j} \wedge dx_{b_j} +\sum^{l}_{j=1} Jdx_{\breve{b}_j} \wedge dx_{\breve{b}_j} + Jdx_{m+1} \wedge dx_{m+1}, \\
\end{eqnarray*}
It follows that $<\beta_{\varphi^a},\beta_{\phi^a}> = <\beta_{\breve{\varphi}^b},\beta_{\breve{\varphi}^b}>= m+1$, 
$<\beta_{\varphi^a},\omega> = 1$ and $<\beta_{\breve{\varphi}^b}, \omega>= -1$. The $\beta_{\varphi^a}$'s all belong to the same orbit under the action of $PU(m,1)$ : let $\nn^{+1}_{\cx H^m}$ be this orbit. The $\beta_{\breve{\varphi}^b}$'s also belong to the same orbit under the action of $PU(m,1)$ and we denote it by $\nn^{-1}_{\cx H^m}$.

\subsection{The positive mass theorem.}

\begin{thm}
Let $(M^m,g,J)$ be an asymptotically complex hyperbolic K\"ahler manifold with even complex dimension
and $\scal_g \geq \scal_{\cx H^m}$. We assume that the cohomology class of $\frac{\Omega}{i\pi}$ is 
integral, providing a line bundle $L$ as above, and that $(M,L)$ defines a spin$^c$ structure. 
Then the formula (\ref{theformula}) defines a (possibly infinite) linear functional $\mu_g$ on $\nn_{\cx H^m}$. 
It is well defined by $g$ up to the natural action of $PU(m,1)$, non-negative on $\nn^{+1}_{\cx H^m} \cup \nn^{-1}_{\cx H^m}$ 
and vanishes if and only if $(M^m,g,J)$ is the complex hyperbolic space.
\end{thm}

Note the assumption about $[\Omega]$ and $L$ are of course satisfied when $M$ is contractible. 

\proof
The proof is nearly the same as in the odd dimensional case. The analytical part is completely similar, cf. \cite{BH}. 
So, for every element $\phi$ of $\kk^c_{\cx H^m}$, one may find an harmonic spinor $\psi$ on $M$ that is asymptotic 
to $\phi$ (harmonic means in the kernel of some Dirac operator, cf. \cite{BH}). We may then proceed to the same computation, 
leading to the same mass integral at infinity, $\mu_g$, defined on the linear span of $Q^c(\kk^c)$ in $\nn$. 
From the computations above we know that $Q^c(\kk^c)$ contains $\nn^{+1}_{\cx H^m} \cup \nn^{-1}_{\cx H^m}$. Using the $PU(m,1)$ invariance, it is easy to see that this is the whole $\nn_{\cx H^m}$. The nonnegativity statement is automatic and the fact 
that $\mu_g$, up to $PU(m,1)$,  depends only on $g$ is like in the odd case. 

We are left to justify the rigidity part. If the mass $\mu_g$ vanishes, we know by construction that every element $\phi$ of $\kk^c_{\cx H^m}$ gives rise to an element $\psi$ of $\kk^c_{g}$ that is asymptotic to $\phi$. In particular, for any element $\beta = Q^c_{\cx H^m}(\phi)$ of $\nn^{\pm 1}_{\cx H^m}$, there is an element $\tilde{\beta} = Q^c_g(\psi)$ of $\nn_{g}$ that is asymptotic to $\beta$. As a consequence, $\nn_{g}$ has maximal dimension : $(M,g)$ is complex hyperbolic. 
\endproof

\appendix

\section*{Appendix : an example.}

To build a non-trivial example of asymptotically complex hyperbolic K\"ahler manifold, we use the symplectic point of view of 
\cite{HS}. We work on $\cx^m$, $m \geq 2$ and let $\rho= e^t$ be the Euclidean distance to the origin. Then $d d^c t = \omega_{FS}$ is the pull-back to $\cx^m \backslash \set{0}$ of the Fubini-Study form on $\cx P^{m-1}$. Every $U(m)$-invariant K\"ahler structure on $\cx^m \backslash \set{0}$ admits a K\"ahler form reading
$$
\omega = \psi'(t) dt \wedge d^ct + \psi(t) \omega_{FS},
$$
where $\psi$ is simply the derivative of a radial K\"ahler potential with respect to the variable $t$. The function $\psi$ and $\psi'$ must of course be positive. In fact $\psi$ can be seen as a moment map for the standard action of $\sph^1$ on $\cx^m \backslash \set{0}$ endowed with the symplectic form $\omega$. If $I$ denotes the range of the function $\psi$, we define the ``momentum profile'' $\Theta$ as the positive function $\Theta := \psi' \circ \psi^{-1}$ on $I = ]\inf \psi, \sup \psi[$. The K\"ahler metric extends smoothly near $0$ if and only if $\inf \psi =0$ and $\Theta$ is smooth near $0$ with $\Theta(0)=0$ and $\Theta'(0)=2$. Setting $x := \psi(t) \in I$, the K\"ahler form is given by
$$
\omega = \omega_{\Theta} = \frac{dx \wedge d^c x}{\Theta(x)} + x \; \omega_{FS}
$$
and its scalar curvature is
$$
s_{\Theta}(x) = \frac{2m(m-1)}{x} - \frac{\partial_{xx} (x^{m-1}\Theta(x))}{x^{m-1}}.
$$
The flat metric corresponds to $\Theta(x)=2x$, $x \in \rl_+$, the complex hyperbolic metric to $\Theta(x)=2x + 2x^2$, $x \in \rl_+$ and the Fubini-Study metric on $\cx P^{m}$ to $\Theta(x)=2x - 2x^2$, $x \in [0,1]$.

Let $\Theta_0(x)=2x^2+2x$ be the momentum profile of the complex hyperbolic metric. We wish to build a $U(m)$-invariant K\"ahler metric with complex hyperbolic asymptotic and scalar curvature bounded from below by the scalar curvature of the complex hyperbolic model. We therefore need to find a positive function $\Theta$ defined on $\rl_+$, such that $\Theta(0)=0$, $\Theta'(0)=2$, $\Theta(x)$ is asymptotic to $2x + 2x^2$ as $x$ goes to infinity and $s_{\Theta} \geq s_{\Theta_0}$. Setting $\Theta(x) = \Theta_0(x) -\alpha(x)$, we are lead to find a non zero function $\alpha$ such that $\Theta_0 -\alpha$ is non-negative, $\alpha(0) = \alpha'(0)=0$, $\alpha=o(x^2)$ as $x$ goes to infinity and $\partial_{xx} (x^{m-1}\alpha(x)) \geq 0$. Such a function $\alpha$ can be obtained by choosing
$$
\alpha(x) := x^{1-m} \int_0^x \int_0^y \chi(z) dz dy,
$$
where $\chi$ is a bump function with support inside $[1,+\infty[$ and with unit integral. One can then observe that $\alpha(x)$ vanishes on $[0,1]$ and is equivalent to $x^{2-m}$ as $x$ goes to infinity, while satisfying the desired differential inequality ; besides, the normalization of $\chi$ ensures $\alpha(x) \leq x$ for every $x \geq 0$, which guarantees $\Theta_0 -\alpha \geq 0$. 
The K\"ahler form is then
$$
\omega_{\Theta}= \omega_{\Theta_0} + \left( \frac{\Theta_0(x)}{\Theta(x)} - 1 \right) \; \frac{dx \wedge d^c x}{\Theta_0(x)}
$$ 
with $\frac{\Theta_0(x)}{\Theta(x)} - 1 \sim \frac{x^{-m}}{2}$. To understand the asymptotic of this, we relate $x$ to the initial variable $t_0$ on the complex hyperbolic space. The properties $\frac{dt_0}{dx}=\frac{1}{\Theta_0(x)}$ and $t_0(x=+\infty)=0$ yield $t_0 \sim -\frac{1}{2x}$. The geodesic distance $r$ to $0$ in the complex hyperbolic model is given by $r = \tanh e^{t_0}$ so that $x \sim \frac{e^{2r}}{4}$ at infinity. So our K\"ahler metric satisfies $\omega = \omega_0 + \oo(e^{-2m \, r})$.

\begin{rem*}
The rate of the falloff to the model at infinity that we obtained is optimal in the radial case. In other words, a function $\alpha$ satisfying the desired properties cannot decay faster at infinity. This is due to the fact that the function $x \mapsto x^{m-1} \alpha(x)$ must be convex on $\rl_+$ with zero first order jet at $0$ ; so either it vanishes identically or it grows at least linearly.
\end{rem*}


\end{document}